\documentclass[12pt]{article}
\usepackage{latexsym}
\usepackage{amsfonts}
\def\C{{\mathbf{C}}}
\def\bC{{\mathbf{\overline{C}}}}

\def\ds{\displaystyle}
\title{Exceptional values in holomorphic families
of entire functions}
\author{Alexandre Eremenko\thanks{
Supported by NSF grants DMS-0100512 and DMS-0244547}} 
\date{March 27, 2005}
\begin{document}
\maketitle

In 1926, Julia \cite{Julia} studied
singularities of implicit functions
defined by equations $f(z,w)=0$, where $f$ is an entire function
of two variables. Among other things, he investigated the exceptional
set $P$ consisting of those $w$ for
which such equation has no solutions $z$.
In other words, $P$ is the complement of the projection
of the analytic set $\{(z,w):f(z,w)=0\}$ onto the second coordinate.
Julia proved that $P$ is closed and
cannot contain a continuum, unless it
coincides with $w$-plane.
Lelong \cite{lelong0} and Tsuji
\cite{tsuji0}, \cite[Thm. VIII.37]{tsuji} independently
improved this result by showing that the logarithmic capacity
of $P$ is zero if $P\neq\C$. 
In the opposite direction, Julia \cite{Julia}
proved that every discrete set $P\subset \C$ can occur as the
exceptional set. 
He writes: {\em ``Resterait \`a voir si cet ensemble, sans
\^etre continu, peut avoir la puissance du continu''}.\footnote{
It remains to see whether this set, without being a continuum, can
have the power of continuum.}

According to Alan Sokal (private communication),
the same question arises
in holomorphic dynamics, when one tries to extend to
holomorphic families of transcendental entire
functions the results of Lyubich \cite{Lyubich}
on holomorphic families of rational
functions.

In this paper, we show that in general,
the result of Lelong and Tsuji is best possible:
every closed set of zero capacity can occur as an exceptional
set (Theorem~1).
Then we study a related problem
of dependence of Picard exceptional values of
the function $z\mapsto f(z,w)$ on the parameter $w$ (Theorem~2).

It is known that the exceptional set $P$
is discrete 
in the
important case that $z\mapsto f(z,w)$ are functions of finite order.
This was discovered by Lelong in \cite{lelong0} and later
the result was
generalized to the case of
multi-dimensional parameter $w$ in \cite[Thm. 3.44]{Lelong}. 

We also mention that the set $P$ has to be analytic in certain
holomorphic families of entire functions with finitely many
singular values, considered in
\cite{EL1,EL2}. These families may consist of functions of
infinite order.

We begin with a simple proof of a version
of Lelong's theorem on functions of finite order.
\vspace{.1in}

\noindent
{\bf Proposition 1.} {\em Let $D$ be a complex manifold,
and $f:\C\times D\to\C$
an analytic function, such that the entire functions
$z\mapsto f(z,w)$ are not identically equal to zero and
are of finite order
for all $w\in D$. Then the set
\begin{equation}
\label{P}
P=\{ w\in D: \,(\forall z\in\C)\, f(z,w)\neq 0\}
\end{equation}
is analytic.}
\vspace{.1in}

\noindent
{\bf Corollary.} {\em Let $D$ be a region in $\C$,
and $f:\C\times D\to\C$ an analytic function,
such that entire functions $z\mapsto f(z,w)$ are
of finite order for all $w\in D$. Then the set $P$ as in
$(\ref{P})$ is discrete or $P=D$.}
\vspace{.1in}

Indeed, the set $A=\{ w\in D: f(.,w)\equiv 0\}$ is discrete
(unless $f=0$ when there is nothing to prove). So there exists
a function $g$ holomorphic in $D$ whose zero set is $A$,
and such that $f/g$ satisfies all conditions of Proposition~1.

\hfill$\Box$
\vspace{.1in}

{\em Remarks.}

1. In general, when $D$ is of dimension greater than one,
and the set $A$ is not empty, one can only prove that
$A\cup P$ is contained in a proper analytic subset of $D$,
unless $A\cup P=D$, \cite[Thm 3.44]{Lelong}.

2. If the order of $f(.,w)$ is finite for all $w\in D$
then this order is bounded on compact subsets of $D$ 
\cite[Thm. 1.41]{Lelong}.
\vspace{.1in}

{\em Proof of Proposition 1.} 
We assume without loss of generality that 
\begin{equation}
\label{one}
f(0,w)\neq 0\quad\mbox{for}\quad w\in D.
\end{equation}
(shift the origin in $\C$ and shrink $D$, if necessary),
and that the order of the function $f(.,w)$ does not
exceed $\lambda$ for all $w\in D$ (see Remark 2 above).

Let $p$ be an integer, $p>\lambda$.
Then, for each $w$, $f$ has the Weierstrass representation
$$\displaystyle
f(z,w)=e^{c_0+\ldots+c_pz^p}\prod_{a:f(a,w)=0}
\left(1-\frac{z}{a}\right)
e^{z/a+\ldots+z^p/pa^p},$$
where $a$ are the zeros of $f(.,w)$ repeated according to their
multiplicities, and $c_j$ and $a$ depend on $w$.
Taking the logarithmic derivative, 
differentiating
it $p$ times, and substituting $z=0$,
we obtain for each $w\in D$:
$$F_p(w)=\left.\frac{d^p}{dz^p}\left(\frac{df}{f\,dz}\right)
\right|_{z=0}=
p!\sum_{a:f(a,w)=0} a^{-p-1}.$$
The series in the right hand side is absolutely convergent because
of our choice of $p$.
The functions $F_p$ are holomorphic in $D$, in view of (\ref{one}).
Clearly $w\in P$ implies $F_p(w)=0$ for all $p>\lambda$.
In the opposite direction, $F_p(w)=0$ for all $p>\lambda$
means that all but finitely many
derivatives with respect to $z$ at $z=0$
of the unction $df/fdz$ meromorphic in $\C$ are equal to zero,
so this meromorphic function is a polynomial, and thus
$f(z,w)=\exp(c_0+\ldots+c_pz^p),$ that is $w\in P$.
\hfill$\Box$
\vspace{.1in}

The following result is due to Lelong and Tsuji but they state it
only for the case $\dim D=1$, and we need a multi-dimensional version
in the proof of Theorem 2 below.
\vspace{.1in}

\noindent
{\bf Proposition 2.} {\em Let $D$ be a connected complex manifold,
and $f:\C\times D\to\C$
an analytic function, such that the entire functions
$z\mapsto f(z,w)$ are not identically equal to zero. 
Then the set
$$P=\{ w\in D: \,(\forall z\in\C)\, f(z,w)\neq 0\}$$
is closed and pluripolar or coincides with $D$.}
\vspace{.1in}

{\em Proof.} Suppose that $P\neq D$.
It is enough to show that every point
$w_0\in D$ has a neighborhood $U$ such that $P\cap U$ is pluripolar.
By shifting the origin in $\C_z$ and shrinking $U$ we may assume
that $f(0,w)\neq 0$ for all $w\in U$. Let $r(w)$ be the smallest of
the moduli of zeros of the entire function $z\mapsto f(z,w)$. If this function
has no zeros, we set $r(w)=+\infty$. We are going to prove that
$\log r$ is a continuous plurisuperharmonic function.

First we prove continuity of $r: 
U \to (0,+\infty]$.
Indeed, suppose that $r(w_0)<\infty$, and let 
$k$ be the number of
zeros of $f(.,w_0)$ on $|z|=r(w_0)$, counting multiplicity.
Let $\epsilon>0$ be so small
that the number of zeros of $f(.,w_0)$ in $|z|\leq r(w_0)+\epsilon$ 
equals $k$. Then 
$$\int_{|z|=r(w_0)\pm\epsilon}\frac{df}{fdz}dz=
\left\{\begin{array}{l} 2\pi i k,\\ 0.\end{array}\right.$$
As the integrals depend on $w$ continuously, we conclude that $r$ is
continuous at $w_0$. Consideration of the case $r(w_0)=+\infty$
is similar.

Now we verify that the restriction of $\log r$ to any complex line
is superharmonic. Let $\zeta\to w(\zeta)$ be the equation of such line,
where $w(0)=w_0$. Let $z_0$ be a zero of $f(.,w_0)$ of the smallest
modulus. We verify the inequality for the averages
of $\log r$ over the circles
$|\zeta|=\delta$, where $\delta$ is small enough.
According to the Weierstrass Preparation theorem, the set 
$Q=\{ (z,\zeta):f(z,w(\zeta))=0\}$ is given in a neighborhood of
$(z_0,0)$ by an
equation of the form
$$(z-z_0)^p+b_{p-1}(\zeta)(z-z_0)^{p-1}+\ldots+b_0(\zeta)=0,$$
where $b_j$ are analytic functions in a neighborhood of $0$, $b_j(0)=0$.
We rewrite this as
\begin{equation}
\label{spa1}
z^p+c_{p-1}(\zeta)z^{p-1}+\ldots+c_0(\zeta)=0,
\end{equation}
where $c_j$ are analytic functions, and 
\begin{equation}
\label{3}
c_0(0)=(-z_0)^p.
\end{equation}
Let $V$ be a punctured disc around $0$; we choose its radius so small
that $c_0(\zeta)\neq 0$ in $V$.
For $\zeta$ in $V$, let $z_i(\zeta),\; i=1,\ldots, p$ be
the branches of the multi-valued function $z(\zeta)$ defined by
equation (\ref{spa1}). Then
$$\displaystyle\begin{array}{l}\ds
\frac{1}{2\pi}\int_{-\pi}^{\pi}\log r(w(\delta e^{i\theta}))d\theta
\leq
\frac{1}{2\pi}\int_{-\pi}^{\pi} \min_{i}\log|z_i(\delta e^{i\theta})|d\theta\\
\\
\ds\leq
\frac{1}{2\pi p}\int_{-\pi}^{\pi}\log\prod_{i}|z_i(\delta e^{i\theta})|d\theta
=
\frac{1}{2\pi p}\int_{-\pi}^{\pi}\log|c_0(\delta e^{i\theta})|d\theta\\
\\ \ds=
\log|z_0|=\log r(w_0),\end{array}$$
where we used (\ref{3}) and harmonicity of $\log|c_0|$ in $V\cup\{0\}$.
This completes the proof of plurisuperharmonicity. As
$P=\{ w:\log r(w)=+\infty\}$ we conclude that $P$ is
pluripolar. \hfill$\Box$
\vspace{.1in}

Our first theorem answers the question
of Julia; it shows that the restriction of finiteness
of order 
cannot be removed in Proposition~1, and that Proposition~2
is best possible, at least when $\dim D=1$.
\vspace{.1in}

\noindent
{\bf Theorem 1.} {\em Let $D\subset\C$ be the unit disc,
and $P$ an arbitrary compact subset
of $D$ of zero capacity. Then there exists a holomorphic function
$f:\C\times D\to\C$, such that for every $w\in P$ the
equation $f(z,w)=0$ has no solutions,
and for each $w\in D\backslash P$ it has infinitely many solutions.}
\vspace{.1in}

It is not clear whether a similar result holds with
multidimensional parameter space $D$ and arbitrary closed
pluripolar set $P\subset D$. 
\vspace{.1in}

{\em Proof.} Let $\phi:D\to D\backslash P$ be a universal covering.
Let $S$ be the set of singular points of $\phi$ on the unit circle.
Then $S$ is a closed set of zero Lebesgue measure.

(We recall a simple proof of this fact.
As a bounded analytic function, $\phi$ has radial limits
almost everywhere. It is easy to see that a point where
the radial limit has absolute value $1$ is not a singular
point of $\phi$. Thus the radial limits exist and belong to $P$
almost everywhere on $S$.
Let $u$ be the ``Evans potential''
of $P$, that is a harmonic function in $D\backslash P$,
continuous in $\overline{D}$
and such that $u(\zeta)=0$ for
$\zeta\in\partial D$ and $u(\zeta)=-\infty$ for $\zeta\in P$.
Such function exists for every compact set $P$ of zero capacity.
Now $v=u\circ \phi$ is a negative harmonic function in the unit disc, whose radial
limits on $S$ are equal to $-\infty$,
thus $|S|=0$ by the classical uniqueness theorem.)

According to a theorem of Fatou, (see, for example, \cite[Ch. VI]{Hoffman}),
for every set $S$ of zero Lebesgue measure on $\partial D$, there exists a holomorphic function
$g$ in $D$, continuous in $\overline{D}$ and such that 
$$\{\zeta\in\overline{D}: g(\zeta)=0\}=S.$$ In particular,
$g$ has no zeros in $D$.  

Now we define the following set $Q\subset\C\times D$:
$$Q=\{(1/g(\zeta),\phi(\zeta)):\zeta\in D\}.$$
It is evident that the projection of $Q$ on the second coordinate equals
$D\backslash P$. It remains to prove that the set $Q$ is analytic.
For this, it is enough to establish that the map
$$\Phi: D\to \C\times D,\quad \Phi(\zeta)=(1/g(\zeta),\phi(\zeta))$$
is proper. Let $K\subset \C\times D$ be a compact subset. Then the closure of
$\Phi^{-1}(K)$ in $\overline{D}$ is disjoint from $S$ because $g$ is continuous and
$1/g(\zeta)\to\infty$ as $\zeta\to S$. On the other hand, for every point $\zeta\in\partial
D\backslash S$, the limit 
$$\lim_{\zeta\to\zeta_0}\phi(\zeta)$$
exists and has absolute value $1$. So $\Phi^{-1}(K)$ is compact in $D$.
 
Now the existence of the required
function $f$ follows from the solvability
of the Second Cousin problem \cite[sect. 5.6]{Horm}.
\hfill$\Box$
\vspace{.1in}

Notice that the map $\Phi$ constructed in the proof is an immersion.
\vspace{.1in}

We recall that a point $a\in\C$ is
called an exceptional value of an entire function $f$
if the equation $f(z)=a$ has no solutions.
Picard's Little theorem says that a non-constant
entire function can have at most one exceptional value.

Let $f$ be an entire function of $z$ depending 
of the parameter $w$ holomorphically, as in Propositions 1 and 2, and
$$\mbox{{\em assume in the
rest of the paper
that for all $w\in D$, 
$f(.,w)\neq{\mathrm{const}}$.}}$$
Let $n(w)\in\{0,1\}$ be the number of exceptional values of $f(.,w)$.
\vspace{.1in}

{\em Question 1.} What can be said about $n(w)$ as a function
of $w$?
\vspace{.1in}

\noindent
{\bf Example 1.} $f(z,w)=e^z+wz$. We have $n(0)=1$ and $n(w)=0$
for $w\neq 0$.
\vspace{.1in}

\noindent
{\bf Example 2.} $f(z,w)=(e^{wz}-1)/w$ for $w\neq 0$, and
$f(z,0)=z$. We have $n(0)=0$,
while $n(w)=1$ for $w\neq 0$.
The exceptional value $a(w)=-1/w$ tends to infinity
as $w\to 0$.
\vspace{.1in}

Thus $n$ is neither upper nor lower semicontinuous.
\vspace{.1in}

\noindent
{\bf Example 3.} Let
$$f(z,w)=\int_{-\infty}^z(\zeta+w)e^{-\zeta^2/2}d\zeta,$$
where the contour of integration consists of the negative ray, passed
left to right,
followed by a curve from $0$ to $z$. We have $f(z,0)=-e^{-z^2/2}$,
which has exceptional value $0$, so $n(0)=1$. It is easy to see
that there are no exceptional values for $w\neq 0$,
so $n(w)=0$ for $w\neq 0$.
Thus $n(w)$ is the same as in Example~1, but this time we have
an additional feature that the set of singular values of
$f(.,w)$ is finite for all $w\in\C$,
namely, there is one critical value
$f(-w,w)$ and two asymptotic values, $0$ and $\sqrt{2\pi}w$.
\vspace{.1in}

{\em Question 2.} Suppose that $n(w)\equiv 1$, and let $a(w)$ be the exceptional
value of $f(.,w)$.
What can be said about $a(w)$ as a function of $w$?
\vspace{.1in}

For functions of finite order, Question 2 was addressed by
Nishino \cite{Nishino} who proved the following:
\vspace{.1in}

\noindent
{\em Let $f$ be an entire function of two variables
such that $z\mapsto f(z,w)$ is a non-constant function
of finite order for all $w$.
If $n(w)=1$ for all $w$ in some set having a finite accumulation
point, then there exists a meromorphic function $\tilde{a}(w)$
such that $a(w)=\tilde{a}(w)$ when $a(w)\neq\infty$,
and  $f(.,w)$ is a polynomial when $\tilde{a}
(w)=\infty.$}
\vspace{.1in}

\noindent
{\bf Example 4.} (Nishino) Let $f(z,w)=(e^{we^z}-1)/w,\; w\neq 0$
and
$f(z,0)=e^z.$ Then $f$ is an entire function of two variables
and $n(w)\equiv 1$. However $a(w)=-1/w,\; w\neq 0$ and $a(0)=0$,
so $a$ is a discontinuous function of $w$.
\vspace{.1in}

Nishino also proved that for arbitrary entire
function $f$ of two variables, with $n(w)=1$ in some region $D$, 
the set of discontinuity of
the function $a(w)$ is closed and nowhere dense in $D$. 

Our Theorem 2 below gives a complete answer to Question 2. 
We first prove the following semi-continuity property of the
set of exceptional values which holds for all
{\em meromorphic functions} holomorphically depending
on parameter. We always assume that $z\mapsto f(z,w)$
is non-constant for all $w$. Denote
$$A(w)=\{ a\in\bC:\, (\forall z\in\C)\, f(z,w)\neq a\}.$$

\noindent
{\bf Proposition 3.} {\em For every $w_0\in D$ and every $\epsilon>0$
there exists $\delta>0$ such that $|w-w_0|<\delta$
implies that $A(w)$ is contained in the $\epsilon$-neighborhood
of $A(w_0)$ with respect to the spherical metric.}
\vspace{.1in}

{\em Proof.} Let $U$ be the open $\epsilon$-neighborhood of $A(w_0)$.
Then $K=\bC\backslash U$ is compact, so there exists $r>0$
such that the image of the disc $|z|<r$ under $f(.,w_0)$ contains $K$.
Then by Hurwitz's theorem, for every $w$ close enough to $w_0$
the image of the disc $|z|<r$ under $f(.,w)$ will also contain $K$.
\hfill$\Box$
\vspace{.1in}

\noindent
{\bf Corollary 1.} {\em The set of meromorphic functions
having no exceptional values
on the Riemann sphere is open in the topology of
uniform convergence on compact subsets of $\C$ with
respect to the spherical metric in the image.}
\vspace{.1in}

The set of entire functions whose only exceptional value is
$\infty$ is not open as Example~2 above shows.
\vspace{.1in}

\noindent
{\bf Corollary 2.} {\em Suppose that $f(.,w)$ is entire
and has an exceptional value
$a(w)\in\C$ for all $w$ on some
subset $E\subset D$. If $a(w)$ is bounded on $E$
then its restriction on $E$ is continuous.}
\vspace{.1in}

Example 4 shows that $a$ can be discontinuous.
\vspace{.1in}

\noindent
{\bf Theorem 2.} {\em Let $f:\C\times D\to\C$,
be a holomorphic function,
where $D$ is a region in $\C$,
and $z\mapsto f(z,w)$ is not constant for
all $w\in D$. Assume that for some function $a:D\to\C$ we have
$f(z,w)\neq a(w)$ for all $z\in\C$.

Then there exists a discrete set $E\subset D$
such that $a$ is holomorphic
in $D\backslash E$, and $a(w)\to\infty$ as $w\to w_0$,
for every $w_0\in E$.}

\vspace{.1in}

So the singularities of $a$ can be only
of the type described in Example 4.

The main ingredient of the proof of Theorem~2
is the following recent
result of N. Shcherbina \cite{Shch}:
\vspace{.1in}

\noindent
{\bf Theorem A.} {\em Let $h$ be a continuous function
in a region $G\in\C^n$. If the graph of $h$ is a pluripolar subset
of $\C^{n+1}$ then $h$ is analytic.}
\vspace{.1in}

{\em Proof of Theorem 2}. Consider the analytic set
$$Q=\{(z,w,a)\in\C_z\times D\times\C_a: f(z,w)-a=0\}.$$
Let $R$ be the complement of the projection of $Q$ onto $D\times\C_a$.
Then it follows from the assumptions of Theorem 2 and Picard's
theorem that $R$
is the graph of the function $w\mapsto a(w)$.
So $R$ is a non-empty proper subset of $D\times\C_a$.
Proposition 1 implies that $R$ is closed in $D$. 

According to Proposition 2,  $R$ is pluripolar.

It follows from Proposition 3 that $w\mapsto |a(w)|$ is lower semi-continuous,
so the sets
$$E_n=\{ w\in D:|a(w)|\leq n\},\quad\mbox{where}\quad n=1,2,3\ldots$$
are closed. We have $E_1\subset E_2\subset\ldots$, and the assumptions
of Theorem 2 imply that $D=\cup E_j$.
Let $E_j^0$ be the interiors of $E_j$, then $E_1^0\subset E_2^0\subset\ldots$,
and the set $G=\cup E_j^0$ is open.

We claim that $G$ is dense in $D$.
Indeed, otherwise there would exist a disc $U\subset D$ which is disjoint from $G$.
But then the closed sets $E_j\backslash E^0_j$ with empty interiors would
cover $D$, which is impossible by the Baire category theorem.
This proves the claim.

As $a$ is locally bounded on $G$, Corollary 2 of Proposition 3
implies that $a$ is continuous
in $G$, so by Shcherbina's theorem, $a$ is analytic in $G$. It follows from the
definition of $G$ that  
$a$ does not have an analytic continuation from any component of $G$ to any boundary
point of $G$. Our goal
is to prove that $D\backslash G$ consists of isolated points.

If $G$ has an isolated boundary point $w_0$
then we have
\begin{equation}
\label{lim}
\lim_{w\to w_0}a(w)=\infty.
\end{equation}
Indeed, by Proposition 3, the limit set of $a(w)$ as $w\to w_0,\; w\in G$
consists of at most two points, $a(w_0)$ and $\infty$.
On the other hand, this limit set is connected.
So the limit exists. If the limit is finite, then it is
equal to $a(w_0)$ and the Removable Singularity Theorem
gives an analytic continuation of $a$ to $G\cup w_0$,
contradicting the statement
above that there is no such continuation.
So the limit is infinite and (\ref{lim}) holds.

We add to $G$ all its isolated boundary points, thus obtaining
new open set $G'$ containing $G$.
Our function $a$ has a {\em meromorphic} continuation
to $G'$, which we call $\tilde{a}$. This meromorphic continuation
coincides with $a$ in $G$ but does not coincide at the added points
$G'\backslash G$.

We claim that $G'$ has no isolated boundary points.
Indeed, suppose that $w_0$ is an isolated boundary point of $G'$.
By the same argument as above, the limit (\ref{lim}) exists and
is infinite. Then $\tilde{a}$ can be extended to $w_0$ such that
the extended function has 
a pole at $w_0$, but then $w_0$ would be an isolated
boundary point of $G$ (poles cannot accumulate to a pole)
so $w_0\in G'$ by definition of $G'$.
This contradiction proves the claim.

Let $F'$ be the complement of $G'$ with respect to $D$.
Then $F'$ is closed, nowhere dense subset of $D$. Furthermore,
$F'$ has no isolated points, because such points would be 
isolated boundary points of $G'$. So $F'$ is perfect or empty.
Our goal is to prove that $F'$ is empty.

Assume the contrary, that is that $F'$ is perfect.
The closed sets $E_n$ cover the locally compact space $F'$, so
by the Baire category theorem one of these $E_n$ contains a
relatively open part of $F'$. This means that there exists a
positive integer $n$ and an open
disc $U\subset D$ intersecting $F'$ and such that
\begin{equation}
\label{n}
|a(w)|\leq n\quad\mbox{for}\quad w\in U\cap F'.
\end{equation}
By Corollary 2 of Proposition 3, this implies that the
restriction of $a$ on $U\cap F'$ is continuous.

We are going to prove that $\tilde{a}$ has a continuous extension
from $G'$ to $U\cap F'$, and this extension
agrees with the restriction of $a$ on $U\cap F'$.
Let $W$ be a point of $U\cap F'$, and $(w_k)$ a sequence in $G'$
tending to $W$.
Choosing a subsequence we may assume that there exists a limit
\begin{equation}\label{limit}
\lim_{k\to\infty}\tilde{a}(w_k),
\end{equation}
finite or infinite.
By a small perturbation of the sequence that does not change
the limit of $\tilde{a}(w_k)$ we may assume that
$w_k$ are not poles of $\tilde{a}$ so $a(w_k)=\tilde{a}(w_k)$.
By Proposition 3, the limit (\ref{limit})
can only be $a(W)$ or $\infty$.

To prove continuity we have to exclude the latter case.
So suppose that
\begin{equation}\label{inft}
\lim_{k\to\infty}\tilde{a}(w_k)=\infty.
\end{equation}
Let $C_k$ be the component of
the set
$$G'\cap\{ w\in U:|w-W|<2|w_k-W|\},$$
that contains $w_k$, and let
\begin{equation}\label{m}
m_k=\inf\{|\tilde{a}(w)|:w\in C_k\}.
\end{equation}
We claim that $m_k\to\infty$. Indeed, suppose this is not so,
then choosing a subsequence we may assume that $m_k\leq m$ for
some $m>n$. Then there exists a curve in $C_k$ connecting $w_k$ to
some point $w_k^\prime\in C_k$ such that $|\tilde{a}(w_k^\prime)|\leq m+1$.
As $|a|$ is continuous in $C_k$, (\ref{inft}) implies that
this curve contains a point $y_k$
such that $|\tilde{a}(y_k)|=m+1.$ By selecting another subsequence
we achieve that
$$\lim_{k\to\infty}\tilde{a}(y_k)=y,
\quad\mbox{where}\quad |y|=m+1>n.$$
As $|a(W)|\leq n$, we obtain a contradiction with Proposition~3.
This proves our claim that $m_k\to\infty$ in (\ref{m}).

So we can assume that 
\begin{equation}
\label{mk}
m_k\geq n+1\quad\mbox{for all}\quad k.
\end{equation}
Let us show that this leads
to a contradiction. Fix $k$, and consider the limit set
of $\tilde{a}(w)$ as $w\to\partial C_k\cap U$ from $C_k\cap U$.
In view of (\ref{n}), (\ref{mk}) and Proposition~3, this limit
set consists of the single point, namely $\infty$.
To see that this is impossible, we use the following
\vspace{.1in}

\noindent
{\bf Lemma.} {\em Let $V$ and $C$ be two intersecting
regions in $\C$,
and $g$ a meromorphic function in $C$ such that
$$\lim_{w\to W,\; w\in C}
g(w)=\infty\quad\mbox{for all}\quad W\in V\cap\partial C.$$
Then $V\cap\partial C$ consists of isolated points in $V$ and $g$ has
meromorphic extension from $C$ to $V\cup C$.}
\vspace{.1in}

{\em Proof}. By shrinking $V$, we may assume that
$|g(w)|\geq 1$ for $w\in C\cap V$
Then $h=1/g$ has a continuous extension from $C$ to $C\cup V$
by setting $h(w)=0$ for $w\in V\backslash C$.
The extended function is holomorphic on the set
$$\{ w\in C\cap V: h(w)\neq 0\},$$
so by Rado's theorem
\cite[Thm. 3.6.5]{Ransford}, $h$ is analytic in $V\cup C$,
$1/g$ gives the required meromorphic extension of $g$.
\hfill$\Box$
\vspace{.1in}

Applying this Lemma with $$C=C_k,
V=\{ w\in U:|w-W|<2|w_k-W|\}$$ and $g=\tilde{a}$,
and taking into account that $\partial C_k\cap V\subset F'$, and
$F'$ has no
isolated points, 
we arrive at a contradiction which completes the proof that
$\tilde{a}$ has a continuous extension to $F'\cap U$, an  extension
which agrees with $a$ on $F'\cap U$.

{}From Theorem A we obtain now that $a$ is analytic on $F'\cap U$,
which contradicts the fact stated in the beginning of
the proof that $a$ has no analytic continuation from $G$.

This contradiction shows that $F'=\emptyset$ which proves the
theorem.
\hfill$\Box$
\vspace{.1in}

The author thanks Alan Sokal for asking the question which
is the subject of this paper, and Adam Epstein, 
Laszlo Lempert, Pietro Poggi-Corradini,
Alexander Rashkovskii and Nikolay Scherbina
for stimulating discussions.

{\em Purdue University, West Lafayette, IN, U.S.A.

eremenko{@}math.purdue.edu}

\end{document}